# TENSOR MODELS, A QUANTUM FIELD THEORETICAL PARTICULARIZATION

Adrian TANASA[1,2]

[1] LIPN, Institut Galilée, CNRS UMR 7030, Université Paris-Nord, 99 av. Clement, 93430 Villetaneuse, France
[2] "Horia Hulubei" National Institute for Physics and Nuclear Engineering, P. O. B. MG-6, 077125 Magurele, Romania
E-mail: adrian.tanasa@ens-lyon.org

Random matrix models have been intensively studied in mathematical physics and have proven useful in combinatorics. In this paper we introduce a generalization of these models to a class of tensor models. As the topology and combinatorics of these tensor models are much more complex than those of matrix models, a first quantum field theoretical (QFT) proposition for simplifying them was recently proposed in R. Gurau (2011). Here, we present a different QFT simplification. We also identify some of its combinatorial properties, in comparison with those of the simplification of R. Gurau (2011).

## 1. INTRODUCTION AND MOTIVATION

The interplay between combinatorics and theoretical physics has increased over the last few decades. This can be justified by the fact that a good knowledge of combinatorics permits theoretical physicists to propose better models, or to better formulate fundamental questions. Thus, combinatorics plays an important rôle in domains of theoretical physics such as statistical physics and integrable systems, but also in general theoretical physical frameworks like that of quantum field theorety (QFT). QFT can be understood as a quantum description of particles and their interactions; this mathematical description is also compatible with Einstein's theory of special relativity. Thus, within the framework of elementary particle physics (or high energy physics), QFT led to the Standard Model, which is the physical theory tested with the best accuracy by experiments. Moreover, the QFT formalism can be successfully applied to statistical physics, condensed matter physics and so on. Combinatorics also plays a non-trivial rôle in such an intricate formalism (see, for example, the thesis [1] for more details).

In this paper, we first present a QFT framework for a potential theory of quantum gravity, the tensor models (see also [2] for a short review of various combinatorial aspects of these tensor models). The main idea behind this approach is to generalize matrix models, models which were proven to be useful in describing (for example) quantum gravity in two dimensions (of space-time). Let us also mention that matrix models are known to be closely related to combinatorial maps, intensively studied in enumerative combinatorics. This was originally motivated by the four colors theorem but it then become an independent subject, with a large number of algorithmic or algebraic applications (see for example the book [3]). The topology and combinatorics of the three-dimensional (or four-dimensional) framework are much more complex than those of the two-dimensional case. A simplification of these models - colored tensor models - has recently been proposed and intensively studied.

In this paper, we present a different proposition [4] for such a simplification of the tensor models, a proposition which follows a QFT logics, as it is explained in the following sections. Several combinatorial properties of this model are presented, in relation to the properties of colored tensor models.

The paper is organized as follows. In the following section, we present some of the main ideas of the QFT formalism. In the third section, we introduce tensor models, as a QFT generalization of matrix models. Moreover, colorable tensor models are defined. The next section presents the main results of this paper, namely the introduction of another QFT particularization and the presentation of some its combinatorial properties (in comparison with those of colorable models). The fifth section gives a possible generalization to the four-dimensional case. The last section presents some concluding remarks and perspectives for future work.

## 2. THE QUANTUM FIELD THEORETICAL FORMALISM IN A NUTSHELL - BUILT-IN COMBINATORICS

We now give a glimpse of some of the aspects of the mathematical formalism of QFT and show why combinatorics is built-in in this formalism. The interested reader may report to various books on the subject (see, for example, [5]). To illustrate these concepts, we choose to work with the simplest QFT model, the $\Phi^4$ one, in which one has a single type of field (that, from a mathematical point of view, is some function)

$$\Phi : \mathsf{R}^D \to \mathsf{K}, \tag{1}$$

where $D \in \mathsf{N}$ and $\mathsf{K}$ is taken to be $\mathsf{R}$ or $\mathsf{C}$. The parameter $D$ is the dimension of the space-time on which the field lives in and is thus taken to be four (three spatial dimensions and one temporal dimension).

A QFT model is defined by means of a functional integral representation of the *partition function* $Z$. Let us now explain in greater deatil what we mean by these notions. One must first defines the *action*, which, from a mathematical point of view, is a functional of the field $\Phi(x)$.

For the $\Phi^4$ model, the action writes:

$$S[\Phi(x)] = \int d^D x \left[ \sum_{\mu=1}^{4} \left( \frac{\partial}{\partial x_\mu} \Phi(x) \right)^2 + \frac{1}{2} m^2 \Phi^2(x) + V[\Phi(x)] \right], \tag{2}$$

where the parameters $m$ and $\lambda$ are the mass and the coupling constant respectively. Moreover, the *interaction potential* writes:

$$V[\Phi(x)] = \frac{\lambda}{4!} \Phi^4(x). \tag{3}$$

It is this term which gives the name, $\Phi^4$, of the respective physical model.

Let us now introduce the functional integration as the product of integrals at each space point $x$ (multiplied by some irrelevant normalization factor): $D\phi(x) := N \prod_x \int d\Phi(x)$. Nevertheless, such an infinite product of Lebesgue measures is mathematically ill-defined. For a well-defined QFT measure, the interested reader may refer, for example, to the review article [6].

The *partition function* is then defined as

$$Z := \int D\Phi(x) e^{-S[\Phi(x)]}. \tag{4}$$

The physical information of a theory is encoded in $n-point\ functions$ (or *correlation functions*) which are defined as:

$$S^{(N)}(x_1, \ldots, x_N) := \frac{1}{Z} \int D\phi(x) \Phi(x_1) \ldots \Phi(x_N) e^{-S[\Phi]} = \langle \Phi(x_1) \ldots \Phi(x_N) \rangle. \tag{5}$$

In general, one is unable to find exact expressions for these correlation functions. In such cases, the tool used in theoretical physics is the *perturbative expansion*, *i. e.* an expansion of the exponential above in powers of $\lambda$. The coefficients of such an expansion are sums of multiple integrals (see (2) and (5)); the number of these integrals grows rapidly with increasing *order* in perturbation theory ( *i. e.* power of the coupling constant $\lambda$ ).

To these multiple integrals, called *Feynman integrals*, one associates *Feynman graphs*, which are very useful in the organization of the expansion coefficients. For the $\Phi^4$ model exhibited here, the graphs have valence four at each vertex; moreover one has internal and external edges (see, for example, Fig. 1).

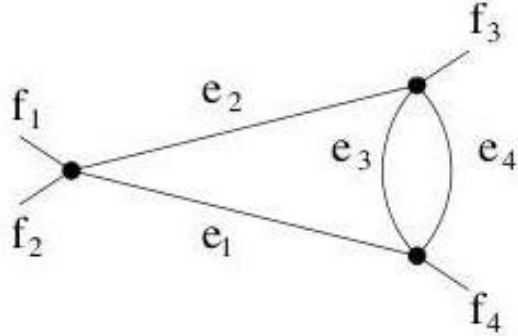

Fig. 1 – A $\Phi^4$ Feynman graph, with four internal edges ($e_1,\ldots,e_4$) and four external edges ($f_1,\ldots,f_4$). It is a graph at the third order in perturbation theory (it has three vertices).

Each of these graphs comes with its combinatorial weight (which is highly non-trivial because of the non-labelling of the edges). The Feynman integrals can then be manipulated using various analytical combinatorial techniques (the Mellin transform or the saddle-point approximation). Let us end this section by giving some insight on *renormalization* in QFT. When computing Feynman integrals, one usually gets different types of divergences (this being a phenomena which appears not only in the $\Phi^4$ model described in this section, but also in more involved QFT models, like gauge theories). The graphs that lead to these divergences then need to be investigated throughly from an analytical point of view. If the respective model is renormalizable, these graphs should correspond to terms present in the action. In order to illustrate what we mean by this, let us get back to the example of the $\Phi^4$ model. The graphs which lead to the various divergences of the model need to have two or four external edges. Thus, these divergences can be ``cured'' by an appropriate renormalization of the parameters of the action (2) (for example, the mass $m$ and the coupling constant $\lambda$). For instance, a graph with four external edges (each of these four external edges being associated to a field $\Phi$) corresponds to the renormalization of the coupling constant $\lambda$, since it is this term in the action which multiplies the $\Phi^4$ term (renormalization of the four-point function, see again (2)). For further details regarding the conditions that a QFT model needs to satisfy in order to be renormalizable, the interested reader can report, for example, to the book [7]. For the sake of completeness, let us also mention that, once the renormalization techniques are performed, the renormalized (and hence finite) correlation functions lead to physical quantities which are measured in elementary particle collider experiments with an extremely high accuracy. For more details on the interplay between combinatorics (topological graph polynomials and combinatorial Connes-Kreimer Hopf algebra of Feynman graphs), the interested reader may report to the survey paper [8].

## 3. TENSOR MODELS; THE COLORABLE MODELS

Let us now go further and introduce tensor models as QFT candidates for a fundamental theory of topological three-dimensional quantum gravity; the interested reader may report to the review paper [9] for further details. The main idea behind this approach is to construct appropriate tensor graphs as duals to triangularization of space-time. In the simplest two-dimensional case, the building block of such a triangularization is of course a triangle. If one takes a dual point of view, the information coming from this triangularization is encoded in some ribbon graph (see, for example, the two-dimensional triangularization of Fig. 2).

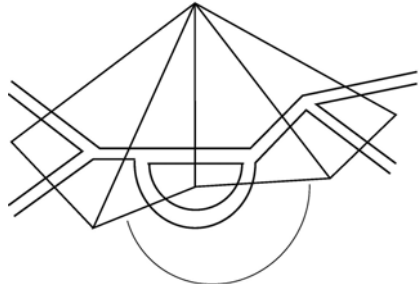

Fig. 2 – Triangularization of a two-dimensional surface and the dual valence three ribbon graph.

The appropriate valence of the associate ribbon graph vertex is thus three: the vertex (see Fig. 3) corresponds to the triangle and the three adjacent edges correspond to the three edges of the original triangle.

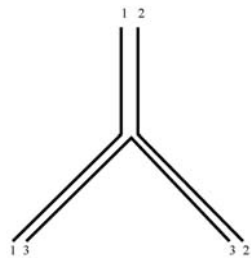

Fig. 3 – A valence three vertex of a two-dimensional quantum gravity model.

The two strands of each edge correspond to the two vertices joined by the respective edge in the initial triangle. The situation is similar for the three-dimensional case. The valence four vertex (see Fig. 4 corresponds to a tetrahedron, the simplest building-block of a three-dimensional space-time.

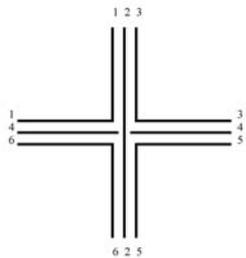

Fig. 4 – A valence four vertex of a three-dimensional tensor model.

The four triangles constructing the tetrahedron correspond to the four edges intersecting at the respective vertex (see, for example, Fig. 5).

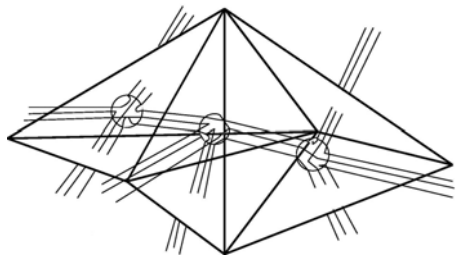

Fig. 5 – Triangularization of a three-dimensional surface and the dual valence four tensor graph.

The three strands of such an edge correspond to the the three edges of the respective triangle (face of the tetrahedron), thus generalizing the image of the two-dimensional case. Finally, the four-dimensional tensor vertex (of valence five), is depicted in Fig. 6.

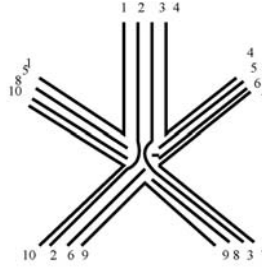

Fig. 6 – A valence five vertex of a four-dimensional quantum gravity tensor model.

These tensor models are thus natural generalizations of the celebrated matrix models, already mentioned in the introduction. The tensor models can thus be associated to three- (and four-)dimensional combinatorial maps, just like matrix models are associated to (the usual two-dimensional) combinatorial maps (the ribbon graphs of matrix models can be seen as a different way of defining maps). Let us now go further with the QFT implementation of these models. We focus on the three-dimensional case, the four-dimensional one being treated in section 5.

As in the $\Phi^4$ model described in the previous section, one has only one type of field, but this time this field does not depend of space-time. One has

$$\phi: G^3 \to \mathbb{R}, \tag{6}$$

where $G$ is some group (these theories are for this reason referred to as group field theory (GFT)). If $G = U(1)$, one speaks of an independent identically distributed (i. i. d.) model, which is the straightforward generalization of the corresponding matrix model. In order to have some connection with the three-dimensional topological quantum gravity, one needs to use as group the holonomy group $SU(2)$ (one can thus refer to GFT as to a theory of holonomies). For the combinatorial issues we deal with in this paper, the choice of the group has no immediate importance. This choice becomes however crucial when computing for example the Feynman integrals associated to these tensor graphs.

The action of this three-dimensional model writes:

$$S[\phi] = \frac{1}{2}\int dg_1 dg_2 dg_3 \phi(g_1,g_2,g_3)\phi(g_1,g_2,g_3) +$$
$$+ \frac{\lambda}{4!}\int dg_1 \ldots dg_6 dg_{1'} \ldots dg_{6'} \phi(g_1,g_2,g_3)\phi(g_{3'},g_4,g_5)\phi(g_{5'},g_{2'},g_6)\phi(g_{6'},g_{4'},g_{1'}) \tag{7}$$
$$\delta(g_1 g_{1'}^{-1})\ldots\delta(g_6 g_{6'}^{-1}),$$

Integrating out the group $\delta-$functions above leads to:

$$S[\phi] = \frac{1}{2}\int dg_1 dg_2 dg_3 \phi(g_1,g_2,g_3)\phi(g_1,g_2,g_3) +$$
$$+ \frac{\lambda}{4!}\int dg_1 \ldots dg_6 \phi(g_1,g_2,g_3)\phi(g_3,g_4,g_5)\phi(g_5,g_2,g_6)\phi(g_6,g_4,g_1). \tag{8}$$

It is this form of the interaction potential which correspond to the graphical representation of Fig. 4. Note that the integration over the group is done with the invariant Haar measure, both in (7) and in (8).

The field $\phi$ in equation (6), taken real-valued, is not assumed here to have specific symmetry properties under the permutations of its arguments: only the identical permutation is associated to such an edge. This type of model is called *orientable* in the GFT literature. The colorable three-dimensional model is defined in the following way. The (real-valued) field (6) is replaced by four complex-valued field $\phi_p$, the index $p = 0,\ldots,3$ being referred to as some color index. Moreover, one now has two types of interactions, a $\phi^4$ one and a $\bar\phi^4$ one. Furthermore, a clockwise cyclic ordering at one of the types of vertices (and an anticlockwise at the second type of vertex) of the four colors at the vertex is also imposed; the action thus writes:

$$S_{\text{col}} = \frac{1}{2}\sum_p \int \bar{\phi}_p \phi_p + \frac{\lambda}{4!}\int \phi_0 \ldots \phi_4 + \frac{\lambda}{4!}\int \bar{\phi}_0 \ldots \bar{\phi}_4. \qquad (9)$$

Note that from now one the integrations over the group are left implicit. Let us also emphasize that the cyclic ordering appearing above has been recently dropped out from this definition being replaced with the constraint that a face of a colorable graph ( *i. e.* a closed circuit) has two colors. Nevertheless, when compairing our results with the one for the colorable models, we will use the initial definition (9).

## 4. ANOTHER QUANTUM FIELD THEORETICAL PARTICULARIZATION; SOME OF ITS COMBINATORIAL FEATURES

Let us now introduce the announced QFT particularization of the model (8). As in the colored case, one has a complex-valued field $\phi$. Nevertheless, we do not need more copies of this field. The proposed interaction is restricted to vertices where each corner has a $+$ or a $-$ label. Furthermore, each vertex has two corners labeled with $+$ and two corners labeled with $-$, which are cyclically ordered as shown in Fig. 7.

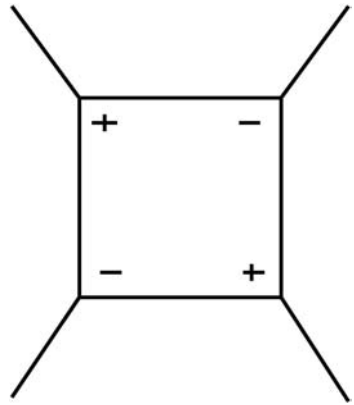

Fig. 7 – The vertex of the proposed model.

A field propagates from a $-$ to a $+$ corner of some vertex.

This notion of corners of a vertex is natural within the framework of non-local QFT Nevertheless, GFT also can be seen in some sense as a non-local QFT on the group manifold on which it lives, since the interaction is not defined on some ``group point'' $g$ (see (7) and (8)), unlike the interaction of the local QFT models (see (2)), where the interaction is $\Phi^4(x)$, localized on some point $x$ of the space-time $\mathsf{R}^4$.

In the original paper [12], this model was called multi-orientable because one has, on the one hand, the orientability of the edges mentioned above (no twists between the various strands) and on the other hand, the orientability of the vertex. Let us also mention that combinatorial maps with four-valence vertices like the orientable ones of Fig. 7 have not been counted in combinatorics [10].

The action of the proposed model writes:

$$S[\phi] = \frac{1}{2}\int \bar{\phi}\phi + \frac{\lambda}{4!}\int \bar{\phi}\phi\bar{\phi}\phi. \qquad (10)$$

As already explained in section 2, a QFT action leads, through an appropriate perturbative expansion, to a certain class of graphs. Trough this mechanism, colorability discards a highly significant class of graphs, including the so-called "wrapping singularities", which correspond to graphs containing loops, *i. e.* edges starting and ending on the same vertex (they are called *tadpoles* in the QFT terminology), like Figs. 8 and 9 (see, for example, [11] and references within for more details on this types of singularities).

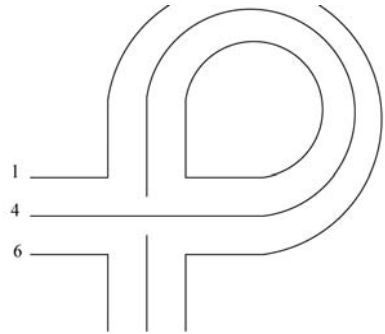

Fig. 8 – A tadpole graph which is generated by the QFT action (10). This graph is also not allowed by the colorable GFT. The indices correspond to the group elements associated to the two external edges.

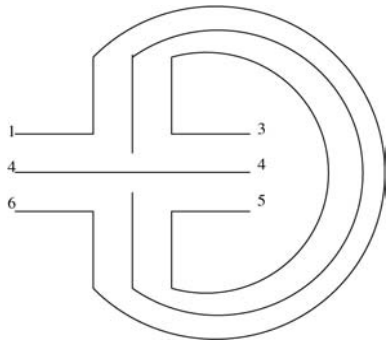

Fig. 9 – A tadpole graph which is not generated by the QFT action (10). The indices correspond to the group elements associated to the two external edges.

**PROPOSITION 4.1** (Proposition $3.1$ of [12]). *Every GFT graph which is colorable is also generated by the QFT action (10).*

The reciprocal statement is not true. A counterexample is the tadpole graph of Fig. 8. Another example, which is not a tadpole graph, is the graph of Fig. 10.

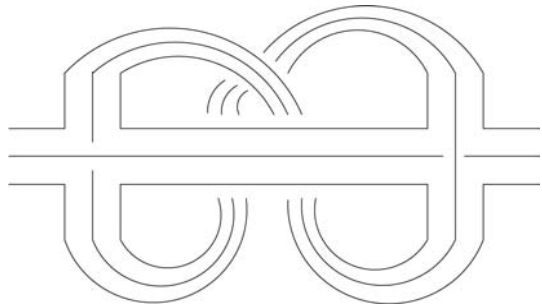

Fig. 10 – A non-tadpole example of a GFT graph.

This can be rephared as: the proposed model (10) discards a less important class of graphs then the colorable model.

### 4.1. Tadpoles and generalized tadpoles

Before investigating the issue of tadpoles within the QFT framework (10), let us make the following remark. Since the edges of the models we deal with here do not allows twists (as already stated above), one can drop the "middle" strand of any edge and obtain an one-to-one correspondence with some ribbon graph (or combinatorial map). The ribbon graph thus obtained is the *jacket* introduced in [12]. One can thus refer to the planarity of the respective tensor GFT graph as to the planarity of the ribbon graph associated in this way. Moreover, one can count the number of faces broken by the external legs; we denote this number by $B$. If $B > 1$ we call the respective graph *irregular*. The tadpole in Fig. 8 is thus referred to as a *planar* tadpole, while the tadpole in Fig. 9 is referred to as a "*non-planar*" tadpole (although the terminology "non-planar" usually used in the QFT literature is inexact, because the respective ribbon graph is planar but it just has a number of broken faces superior to one). We have seen above that planar tadpoles are allowed by the model (10), while the non-planar ones are not.

Let us now recall from [13] the following definition:

**Definition 4.1.** *A generalized tadpole is a graph with one external vertex.*

Planar generalized tadpoles are allowed by models like the one defined in (10). A ``non-planar'' generalized tadpole is a graph with two external edges and with $B = 2$. These graphs are not allowed by models of type (10).

### 4.2. Tadfaces

We recall from [12] the following definition:

**Definition 4.2.** *A tadface is a face that goes at least twice through a line.*

Let us give a few more explanation of this. Such a tadface can be obtained if, one goes through the respective edge a first time through a strand and a second time through a second strand. A second-order (in perturbation theory) example of such a graph is given in Fig. 11. Nevertheless, this graph is not allowed by the model (10), since it is made up of two "non-planar" tadpoles like the ones in Fig. 9 (each of them not being allowed by the model (10), as already explained above).

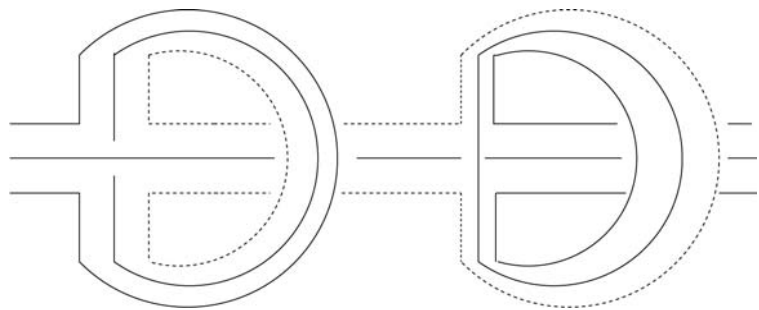

Fig. 11 – A tadface graph. The tadface is represented by the dashed line. One notices that the respective face goes twice through the edge relying the two tadpoles, once through one strand and once through another strand of the edge.

One has the following result:

**THEOREM 4.1** (Theorem $3.1$ of [4]). *Tadfaces are not allowed by the QFT model* (10).

One can also obtain the following result, already announced in [11]:

**COROLLARY 4.1** (Corollary $3.1$ of [4]). *Tadfaces are not allowed by colorability*.

Furthermore, one has:

**COROLLARY 4.2** (Corollary $3.2$ of [4]). *"Non-planar" generalized tadpoles are not allowed by multi-orientability*.

Finally, one has:

**COROLLARY 4.3** (Corollary $3.3$ of [4]). *"Non-planar" generalized tadpoles are not allowed by colorability*.

We resume all these results in the following table, which compares the two models

|  | colorable | model (10) |
| --- | --- | --- |
| generalized planar tadpoles | forbidden | allowed |
| generalized ``non-planar'' tadpoles | forbidden | forbidden |
| tadfaces | forbidden | forbidden |

For the sake of completeness, let us mention that in the original paper [4], some Feynman integral computations have been performed. Furthermore, the comparison with the colorable model was pushed further by investigated aspects related to the renormalizability (see section 2) of the two types of GFT models. Thus, it was exhibited that there exists graphs which, within the colorable framework, are divergent (when one computes the associated Feynman integral) but they do not correspond to terms present in the action (9). For the QFT framework defined by the model (10), this is cured, in the sense that these graphs, which are still divergent, correspond however to terms present in the action (10) (see again the original paper [4] for more details). An example of such a graph is given in Fig. 12.

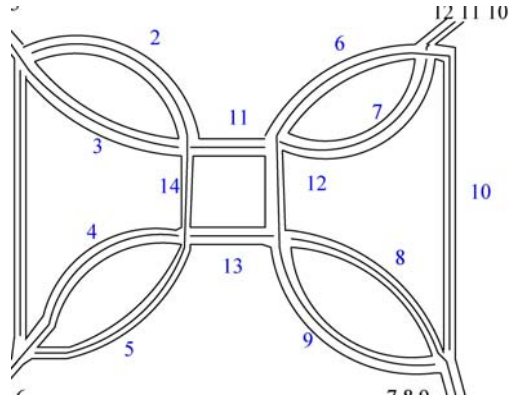

Fig. 12 – An example of a tensor graph which represents a contribution of a form not present in the bare action. In the QFT framework of the model (10), it represent a contribution of a form already present in the bare action. The fourteen internal edges and the the twelve group elements associated to the four external edges have been labeled.

## 5. GENERALIZATION TO THE FOUR-DIMENSIONAL CASE

The generalization of the model (10) to GFT in even dimension is not straightforward. This comes from the fact that the interaction $\phi^{D+1}$ is odd; one thus has two inequivalent choices of distribution of the $-$ and $+$ signs on the corners of the vertex. For the four-dimensional case (where the vertex is given in Fig. 6), case of interest for quantum gravity, these two possibilities of interactions are given in Fig. 13 and 14.

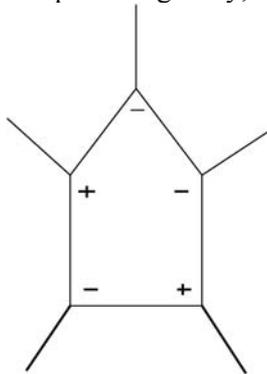

Fig. 13 – A first possibility of a vertex for four-dimensional GFT.

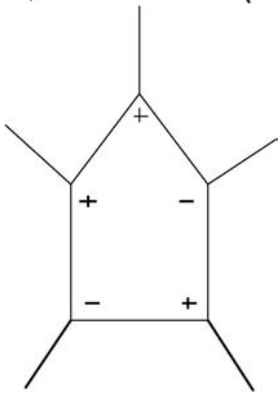

Fig. 14 – A second possibility of a vertex for four-dimensional GFT.

The action of the proposed model thus writes

$$S[\phi] = \frac{1}{2}\int \bar{\phi}\phi + \frac{\lambda_1}{5!}\int \bar{\phi}\phi\bar{\phi}\bar{\phi}\bar{\phi} + \frac{\lambda_2}{5!}\int \bar{\phi}\phi\bar{\phi}\bar{\phi}\phi, \quad (11)$$

where, as in (9) or (10), the integrations over the group are left implicit. Moreover, in order to keep generality, we choose to have two distinct coupling constants $\lambda_1$ and $\lambda_2$. The considerations of the previous section extend for this four-dimensional model.

## 5. CONCLUSION AND PERSPECTIVES

We have introduced in this paper the model (10), as a way of simplifying the topology and combinatorics of GFT; this simplification is QFT-inspired and is different of the one proposed within the colorability of GFT. An analysis of the differences between the classes of tensor graphs discarded by these two types of models has been done. A generalization from three-dimensional to four-dimensional models has also been proposed. Since this is a proposal for a new type of GFT models, the perspectives for future work on this subject appears to us as particularly vast. Thus, it would be interesting to check whether or not the various achievements obtained within the framework of colorable GFT models (see, for example, [11] and references within) can be adapted (and in what conditions) to the model (10). Among these perspectives, we can thus list here the definition of some computable cellular homology, as well as the identification of the class of dominant tensor graphs (in the two-dimensional case, this rôle is played by the graphs corresponding to triangularization of the sphere).

Finally, let us mention here as a related perspective for future work the generalization of the matrix integral techniques, techniques which proved useful, in the two-dimensional case, for map counting (see, for example, [14] for an interesting review).

## ACKNOWLEDGEMENTS

The author acknowledges the CNRS PEPS grant "CombGraph". The grants PN 09 37 01 02 and CNCSIS Tinere Echipe 77/04.08.2010 are also gratefully acknowledged.

## REFERENCES

1. A. TANASA, *Combinatorics in Quantum Field Theory and Random Tensor Models*, 2012, Habilitation thesis, in preparation.
2. A. TANASA, *Combinatorics of random tensor models*, Proceedings of the Romanian Academy A, **13**, pp. 27–31, 2012.
3. S. K. LANDO, A. K. ZVONKIN, *Graphs on surfaces and their application*, Springer, 2004.
4. A. TANASA, *Multi-orientable Group Field Theory*, J. Phys. A, **45**, 165401, 2012; arXiv:1109.0694.
5. H. KLEINERT, V. SCHULTE-FROHLINDE, *Critical Properties of $\Phi^4$ – Theories*, World Scientific, 2000.
6. V. RIVASSEAU, *An introduction to renormalization*, Poincaré Seminar, 2002, pp. 139–177.
7. V. RIVASSEAU, *From perturbative to constructive renormalization*, Princeton University Press, 1992.
8. A. TANASA, *Some combinatorial aspects of quantum field theory*, 2011; arXiv:1102.4231.
9. L. FREIDEL, *Group field theory: An Overview*, Int. J. Theor. Phys., **44**, pp. 1769–1783, 2005.
10. E. FUSY, Private Communication, 2011.
11. R. GURAU, J. RYAN, *Colored tensor models – a review*; arXiv:1109.4812 [hep-th], SIGMA (in press).
12. J. BEN GELOUN, T. KRAJEWSKI, J. MAGNEN, V. RIVASSEAU, *Linearized Group Field Theory and Power Counting Theorems*, Class. Quant. Grav., **27**, 155012, 2010.
13. J. BEN GELOUN, J. MAGNEN, V. RIVASSEAU, *Bosonic Colored Group Field Theory*, Eur. Phys. J. C, **70**, pp. 1119–1130, 2010.
14. A. K. ZVONKIN, *Matrix integrals and map enumeration: An accessible introduction*, Computers and Mathematics with Applications: Mathematical and Computer Modeling, **26**, pp. 281–304, 1997.